\documentclass{elsart}
\usepackage[sumlimits]{amsmath}
\usepackage{amsfonts}
\usepackage{graphicx,graphics}
\DeclareMathOperator{\diff}{d}
\newcommand{\pp}[2]{\frac{\partial #1}{\partial #2}} 

\newcommand{\dd}[2]{\frac{\diff#1}{\diff#2}}

\newcommand{\cmc}{\sum_{i}C^T_i(M^u)^{-1}C_i}
\newcommand{\U}{{\mathrm{u}}}
\newcommand{\h}{{\mathrm{h}}}
\newcommand{\p}{{\mathrm{p}}}
\newcommand{\f}{{\mathrm{f}}}
\newcommand{\F}{{\mathrm{F}}}
\newcommand{\g}{{\mathrm{g}}}
\newcommand{\N}{{\mathrm{N}}}
\usepackage{xspace}
\newcommand{\pdgp}{$\textrm{P1}_{\textrm{DG}}$-$\textrm{P2}$\xspace}
\newcommand{\pdg}{$\textrm{P1}_{\textrm{DG}}$\xspace}
\begin{document}
\begin{frontmatter}
\title{LBB Stability of a Mixed Discontinuous/Continuous Galerkin Finite Element Pair}
\author[aeroIC]{Colin J. Cotter\corauthref{cor}}
\corauth[cor]{Corresponding author.}
\ead{colin.cotter@imperial.ac.uk}
\author[eseIC]{David A. Ham}
\ead{d.ham@imperial.ac.uk}
\author[eseIC]{Christopher C. Pain}
\ead{c.pain@imperial.ac.uk}
\author[potsdam]{Sebastian Reich}
\ead{sreich@math.uni-potsdam.de}
\address[aeroIC]{ Department of Aeronautics,
Imperial College London, London SW7 2AZ, United Kingdom}
\address[eseIC]{Department of Earth Science and Engineering,
Imperial College London, London SW7 2AZ, United Kingdom}
\address[potsdam]{Institut f\"ur
Mathematik, Universit\"at Potsdam, Am Neuen Palais 10, D-14469,
Potsdam, Germany}

\begin{abstract}
  We introduce a new mixed discontinuous/continuous Galerkin finite
  element for solving the 2- and 3-dimensional wave equations and
  equations of incompressible flow. The element, which we refer to as
  \pdgp, uses discontinuous piecewise linear functions for velocity
  and continuous piecewise quadratic functions for pressure. The aim
  of introducing the mixed formulation is to produce a new flexible
  element choice for triangular and tetrahedral meshes which satisfies
  the LBB stability condition and hence has no spurious zero-energy
  modes.  We illustrate this property with numerical integrations of
  the wave equation in two dimensions, an analysis of the resultant
  discrete Laplace operator in two and three dimensions, and a normal
  mode analysis of the semi-discrete wave equation in one dimension.
\end{abstract}

\end{frontmatter}

\section{Introduction}

One of the key strengths of the finite element method is the extensive
choice of element types; this strength leads to endless discussion
amongst practioners about the various benefits of different options.
Alongside issues such as accuracy and practicality, a key issue is
that of LBB stability. This issue manifests itself in the
discretisation of the wave equation (and nonlinear extensions such as
the shallow-water equations and the compressible Euler equations), and
also features in the discretisation of the equations of incompressible
flow. If one considers the wave equation written as a two-component
system
\begin{equation}
\label{wave}
\vec{u}_t + \nabla h = 0, \quad h_t + \nabla\cdot\vec{u} = 0,
\quad \vec{u} = (u_1,\ldots,u_d),
\end{equation}
then finite element discretisation results in 
\[
\dd{}{t}M^u\U_i = -C_i\h , \quad i=1,\ldots d, \quad \dd{}{t}M^h\h
= \sum_{i=1}^dC^T_i\U_i, 
\]
where $C_i$, $i=1,\ldots,d$ are the finite element approximations of
the Cartesian components of the gradient operator,
$-\sum_{i=1}^dC^T_i$ is the finite element approximation to the
divergence operator, $M^u$ and $M^h$ are the mass matrices associated
with the finite element spaces for $\U$ and $\h$ respectively, and $d$
is the number of physical dimensions.  By eliminating $\U$, we obtain
the discrete wave equation
\[
M^h\dd{^2}{t^2}\h - \cmc\h=0.
\]
If the discrete Laplace operator $(M^h)^{-1}\cmc$ has null space of dimension
greater than one, this results in spurious zero-energy solutions which
pollute the solution after a period of time. 

The null space problem also manifests itself in incompressible flow
where the equations consist of a dynamical equation for $\vec{u}$ plus
an incompressibility constraint which is maintained by a pressure
gradient:
\[
\vec{u}_t + N(\vec{u}) = -\nabla p + \vec{F}, \quad \nabla\cdot\vec{u}=0,
\]
where $N$ is the advective nonlinearity and $\vec{F}$ represents all other
forces. In this case the spatial discretisation becomes
\[
M^u\dd{}{t}\U_i + \N_i(\U) = -C_i\p + 
\F_i, \quad i=1,\ldots,d,
\quad \sum_{i=1}^dC^T_i\U_i=0.
\]
The pressure can be obtained by applying $\sum_iC^T_i(M^u)^{-1}$ to
the dynamical equation for $\vec{u}$ to obtain
\[
0 = \dd{}{t}\sum_iC_i^T\U_i = -\cmc\p
-\sum_iC^T_i(M^u)^{-1}(\F_i-\N_i(\U)), \quad i=1,\ldots,d,
\]
which can be solved for $\p$ (after fixing the constant
component $\p_0$) provided that $\cmc$ has a 1-dimensional null
space containing only constant functions.

The analysis of the stability properties of finite element
discretisations associated with spurious eigenvectors of $\cmc$ was
performed by Ladyzhenskaya \cite{La1969}, Babuska \cite{Ba1971} and
Brezzi \cite{Br1974}; as a result, a finite element discretisation is
said to be LBB stable if $\cmc$ is free from spurious eigenvectors. In
general these spurious eigenvectors consist of extra null vectors as
well as ``pesky modes'' which have eigenvalues which tend to zero as
the mesh size goes to zero. The spurious null vectors, which generally
only occur in discretisations using structured meshes, make it
impossible to invert the discrete Laplace matrix. The ``pesky modes'',
which arise on unstructured meshes, are nearly as problematic as they
lead to very large condition numbers for the discrete Laplace matrix
which make iterative methods very slow to converge. For the wave
equation, a finite element discretisation is LBB stable only if the
number of degrees of freedom (DOF) for $h$ is less than the number of
DOF for each component of $\vec{u}$; similarly this specifies a
stability condition on the number of DOF for pressure in
incompressible flow. However, this is not a sufficient condition and
some element choices can still be unstable despite having less DOF for
$h$ than for each component of $\vec{u}$. This stability condition
often leads to staggered grids (such as the $C$-grid in finite
difference/volume methods) and mixed finite elements (see
\cite{Jo2003} for a discussion of mixed elements applied to the wave
equation). The number of DOF for $\vec{u}$ is often increased by
introducing interior modes (``bubble'' functions).  In this paper we
suggest an alternative way of increasing the DOF for $\vec{u}$ by
admitting discontinuous functions (see \cite{KaSh2005} for an review
of discontinuous finite elements and their application to
computational fluid dynamics, see \cite{AiMoMu2006,GrScSc2006} for
applications of discontinuous Galerkin methods to the wave equation
and \cite{Levin2006} for applications to ocean modelling), whilst
keeping the continuity constraint for $h$. This often means that it is
possible to increase the order of accuracy of the discretisation of
$h$ whilst keeping the mixed element LBB stable. For a full treatment
of LBB stability and a summary of the stability properties of a wide
range of element pairs see \cite{GrSa2000}; for an analysis of element
pairs applied to the linearised shallow-water equations see
\cite{Ro_etal2005}.

In section \ref{mixed} we introduce the mixed discontinuous/continuous
\pdgp element in one, two and three dimensions and show how the
boundary conditions are implemented. We also give some values for the
$h$ and $\vec{u}$ DOF which show the effects of making $\vec{u}$
discontinuous. In section \ref{1d} we compute the numerical dispersion
relation for this element applied to the semi-discrete wave equation
which shows that the element is indeed stable in one dimension. The
numerical dispersion relation has a gap in the spectrum between two
branches and we show that the modes from the lower frequency branch
have smaller discontinuities in $\vec{u}$ with the lowest frequencies being
nearly continuous. In section \ref{laplace} we show eigenvalues of
discrete Laplace matrices constructed on various unstructured grids in
two and three dimensions which show that the element is stable. In
section \ref{wave eqn} we show the results of a wave equation calculation
in two dimensions on an unstructured grid which illustrates the
absence of spurious modes. Finally, in section \ref{summary} we give a
summary of the paper and discuss other aspects of this element which
may make it a good choice for ocean modelling applications.

\section{The mixed element}
\label{mixed}
In this section we describe our mixed element formulation in one, two
and three dimensions. 
\subsection{Weak formulation}
We start with the wave equation in the form
(\ref{wave}) with boundary conditions
\begin{equation}
\label{bcs}
\pp{h}{n} = f \quad\mathrm{on}\quad \partial\Omega^N, \quad
h = g \quad\mathrm{on}\quad \partial\Omega^D,
\end{equation}
where $\partial\Omega^N$ and $\partial\Omega^D$ form a partition of
the boundary $\partial\Omega$ of the domain $\Omega$, and multiply by
test functions $\vec{w}$ and $\phi$ and integrate over the whole
of $\Omega$ to obtain
\begin{eqnarray}
\label{weak 1}
\dd{}{t}\int_{\Omega} \vec{w}\cdot\vec{u}\diff{V} & = &
-\int_{\Omega}\vec{w}\cdot\nabla h\diff{V}, \\
\dd{}{t}\int_{\Omega} \phi h \diff{V} & = & 
-\int_{\Omega} \phi\nabla\cdot\vec{u} \diff{V}.
\label{weak 2}
\end{eqnarray}
We then integrate equations (\ref{weak 1})-(\ref{weak 2}) by parts, make
use of the boundary conditions (\ref{bcs}), and finally integrate equation
(\ref{weak 1}) by parts again to obtain
\begin{eqnarray}
\dd{}{t}\int_{\Omega} \vec{w}\cdot\vec{u}\diff{V} & = &
-\int_{\Omega}\vec{w}\cdot\nabla h\diff{V} +
\int_{\partial\Omega^D} \vec{w}\cdot\vec{n} h\diff{S} 
- \int_{\partial\Omega^D} \vec{w}\cdot\vec{n} g \diff{S}, \\
\dd{}{t}\int_{\Omega} \phi h \diff{V} & = & 
\int_{\Omega} \nabla\phi\cdot\vec{u} \diff{V}
- \int_{\partial\Omega^D}\vec{n}\cdot\vec{u}\phi\diff{S}
- \int_{\partial\Omega^N}\phi f\diff{S},
\end{eqnarray}
which is the weak form that we discretise. The key feature of this form
is that derivatives are only applied to the scalar functions $h$ and $\phi$
and not the vector functions $\vec{u}$ and $\vec{w}$ which we shall
discretise with discontinuous elements.
\subsection{The \pdgp element}
\begin{figure}
\begin{center}
\includegraphics*[width=8cm]{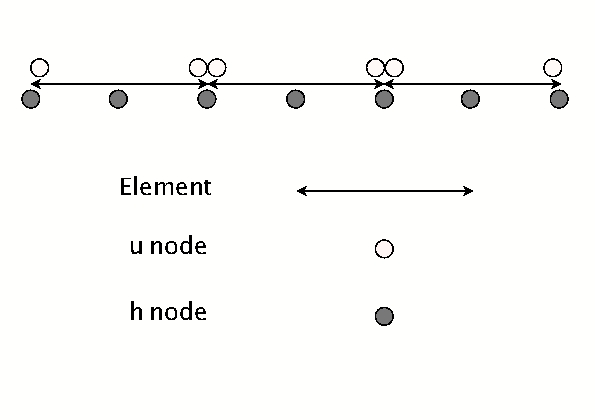}
\includegraphics*[width=4cm]{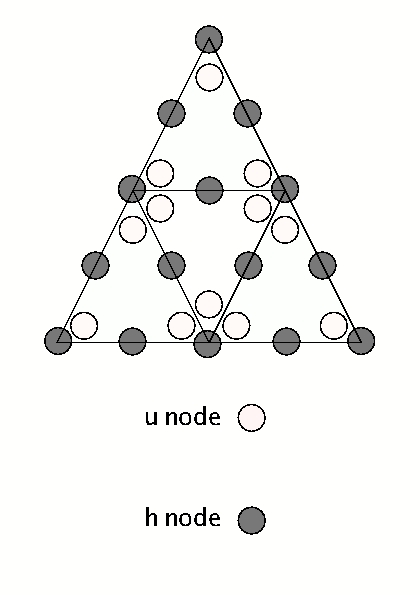}
\end{center}
\caption{\label{1d-dof} Figure showing the DOF for the one-dimensional
  \pdgp element (left) and the two-dimensional \pdgp element
  (right). In one dimension, each element contains two local $\vec{u}$
  DOF and three local $h$ DOF, but the global $h$ DOF are constrained
  to be continuous across element boundaries. In two dimensions there
  are three local $\vec{u}$ DOF and six local $h$ DOF.}
\end{figure}
We then make the choice that the vector functions be discretised with
discontinuous piecewise-linear (\pdg) elements and the scalar
functions be discretised with continuous quadratic (P2) elements. The
reason for this choice is the favourable balance between $\vec{u}$ and
$h$ DOF in this element. 

We write the global finite element expansions in the form
\[
u_i(\vec{x}) = \sum_{\alpha=1}^{m_u}\mathrm{u}_{\alpha,i}N_\alpha(\vec{x}),
\quad h(\vec{x}) = \sum_{\beta=1}^{m_h}\mathrm{h}_{\beta}\bar{N}_\beta(\vec{x}),
\]
where
\[
\U_{1,i} = [\mathrm{u}_{1,i},\ldots,\mathrm{u}_{m_u,i}], \quad i=1,\ldots,d, 
\quad \h = [\mathrm{h}_1,\ldots,\mathrm{h}_{m_h}],
\]
and where $m_u,m_h$ are the numbers of DOF for each component
of $\vec{u}$ and for $h$ respectively. This leads to the following equations:
\[
\dd{}{t}M^u\U_i = -C_i\h - \g_i, \quad i=1,\ldots d, \quad \dd{}{t}M^h\h
= \sum_{i=1}^dC^T_i\U_i - \f,
\]
where
\begin{eqnarray*}
M^u_{\alpha,\beta} & = & \int_{\Omega} N_{\alpha}(\vec{x})N_{\beta}(\vec{x})
\diff{V}(\vec{x}), \\
C_{\alpha,\beta,i} & = & \int_{\Omega} N_{\alpha}(\vec{x})\pp{}{x_i}
\bar{N}_{\alpha}(\vec{x})\diff{V}(\vec{x}) - \int_{\partial\Omega^D}
n_iN_{\alpha}(\vec{x})\bar{N}_{\beta}(\vec{x})\diff{S}, \\
M^h_{\alpha,\beta} & = & \int_{\Omega} \bar{N}_{\alpha}(\vec{x})
\bar{N}_{\beta}(\vec{x})\diff{V}(\vec{x}), \\
g_{\alpha,i} & = & \int_{\partial\Omega^D} N_\alpha(\vec{x}) g(\vec{x})n_i(\vec{x}) \diff{S}, \\
f_{\beta} & = & \int_{\partial\Omega^N} \bar{N}_\beta(\vec{x}) f(\vec{x}) \diff{S},
\end{eqnarray*}
and $d$ is the number of physical dimensions. 

One of the advantages of this element choice is that the mass matrix
$M^u$ is block diagonal (since $\vec{u}$ is discontinuous and so each
global basis function is supported on only one element). This means
that the discrete Laplacian $\cmc$ is still sparse and it is not
necessary to lump the mass matrix when solving the pressure equation
for incompressible flow.

\subsubsection{One dimension} In one dimension on a bounded interval of $I$
elements, there are two local DOF per element for $u$, and so there
are $2I$ global DOF as $\vec{u}$ is discontinuous. There are three
local DOF per element for $h$ but there are $I-1$ global continuity
constraints on the interfaces between each element (see figure
\ref{1d-dof}). This means that there are $3I - (I-1) = 2I + 1$ global
DOF for $h$, and so there is always one more $h$ DOF than $u$
DOF. However, for strong Dirichlet conditions for $h$, or periodic
boundary conditions, we decrease the number of $h$ DOF and gain the
potential for a stable element.

\subsubsection{Two dimensions} In two dimensions we have $F$
triangular elements, with three local DOF per element $\vec{u}$ and
six local DOF per element for $h$. There are no continuity constraints
for $\vec{u}$ and so there are $3F$ DOF (see figure
\ref{1d-dof}). There is an $h$ DOF situated at each vertex and an $h$
DOF situated on each edge, and so there are $V+E$ $h$ global DOF,
where $V$ is the number of vertices and $E$ is the number of
edges. Euler's formula gives $E=V+F+1$ and so there are $2V + F + 1$
$h$ DOF. This means that it is always possible to modify a mesh so
that there are more $\vec{u}$ DOF than $h$ DOF \emph{e.g.} by
repeatedly inserting new vertices into a triangles, breaking them into
four, each time increasing $V$ by 1 and $F$ by 3. In practise, useful
meshes generally satisfy $F>V$ and so this property is
satisfied. Strong Dirichlet boundary conditions for $\vec{u}$ may
reduce the number of $\vec{u}$ DOF below that of $h$. Table \ref{2D
DOF} gives some DOF for various unstructured Delaunay meshes in a
square domain.
\begin{table}
\begin{center}
\begin{tabular}{|c|c|c|c|c|c|}
\hline
Mesh triangles & 36 & 79 & 151 & 1586 & 15574 \\
\hline
Mesh vertices & 24 & 48 & 87 & 820 & 7890 \\
\hline
$\vec{u}$ DOF & 108 & 237 & 453 & 4758 & 46722 \\
\hline
$h$ DOF & 85 & 176 & 326 & 2414 & 31354 \\
\hline
\end{tabular}
\end{center}
\caption{\label{2D DOF} Table showing degrees of freedom for the \pdgp
element pair in two dimensions.  The ratio of $\vec{u}$ DOF to $h$ DOF
appears to converge to 1.5 for large unstructured meshes.}
\end{table}

\subsubsection{Three dimensions} In three dimensions there are four local
$\vec{u}$ DOF and six local $h$ DOF. As there are no continuity
constraints for $\vec{u}$, there are $4T$ global $\vec{u}$ DOF, where
$T$ is the total number of tetrahedra. There is one global $h$ DOF on
each vertex, and one global $h$ DOF on each edge, so there are $V+E$
global $h$ DOF. As for three dimensions, it always possible to
increase $T$ relative to $V+E$ by splitting elements. Table \ref{3D
DOF} gives some DOF for sample unstructured Delaunay meshes in a cubic
domain.
\begin{table}
\begin{center}
\begin{tabular}{|c|c|c|c|c|c|c|}
\hline
Mesh tetrahedra & 44 & 215 & 398 & 2003 & 19140 \\
\hline
Mesh vertices & 26 & 80 & 130 & 488 & 3690 \\
\hline
Mesh edges & 93 & 227 & 633 & 2792 & 24165 \\
\hline
$\vec{u}$ DOF & 176 & 860 & 1592 & 8012 & 77640 \\
\hline
$h$ DOF & 119 & 307 & 763 & 3280 & 27855 \\
\hline
\end{tabular}
\end{center}
\caption{\label{3D DOF} Table showing degrees of freedom for the \pdgp
element pair in three dimensions. The ratio of $\vec{u}$ DOF to $h$ DOF
appears to converge to 2.5 for large unstructured meshes.}
\end{table}

\section{One-dimensional analysis}
\label{1d}
In this section we analyse the \pdgp element applied to the scalar
wave equation in one-dimension on a regular grid with periodic
boundary conditions. 

The local (elemental) mass and gradient matrices are:
\[
\bar{M}_{ij}^u = \int_0^{\Delta x}N_iN_j\diff{x},
\quad \bar{M}_{ij}^h = \int_0^{\Delta x}\bar{N}_i\bar{N}_j\diff{x},
\quad \bar{C}_{ij} = -\int_0^{\Delta x}N_i\dd{}{x}\bar{N}_j\diff{x},
\]
where $\{N_1,N_2\}$ are the linear Lagrange polynomials used to
represent $u$ in the element, $\{N_1,N_2,N_3\}$ are the quadratic
Lagrange polynomials used to represent $h$, $\bar{M}^u$ is the local
mass matrix for $u$, $\bar{M}^h$ is the local mass matrix for $h$ and
$\bar{C}$ is the local gradient matrix. These matrices are
\begin{eqnarray*}
\bar{C} &=& 
\begin{pmatrix}
-5/6 & 2/3 & 1/6 \\
-1/6 & -2/3 & 5/6 \\
\end{pmatrix}, \\
\bar{M}^u &=& \Delta x\begin{pmatrix}
1/3 & 1/6 \\
1/6 & 1/3 \\
\end{pmatrix}, \\
\bar{M}^h &=& \Delta x\begin{pmatrix}
2/15 & 1/15 & -1/30 \\
1/15 & 8/15 & 1/15 \\
-1/30 & 1/15 & 2/15 \\
\end{pmatrix}.
\end{eqnarray*}
After assembling the equations on a regular grid with element width
$\Delta x$, we obtain
\small
\begin{eqnarray}
\label{semi u+}
\frac{\Delta x}{6}\dd{}{t}(2u^n_+ + u^{n+1}_-) & = & - \frac{1}{6}
(-5h^n + 4h^{n+1/2} + h^{n+1}), \\
\frac{\Delta x}{6}\dd{}{t}(u^n_+ + 2u^{n+1}_-) & = & - \frac{1}{6}
(-h^n - 4h^{n+1/2} + 5h^{n+1}), \\
\frac{\Delta x}{30}\dd{}{t}
(-h^{n-1} + 2h^{n-1/2} + 
\qquad & & \nonumber \\
8h^n+ 2h^{n+1/2} - h^{n+1}) & = & 
\frac{1}{6}(u^{n-1}_+ + 5u^n_--5u^n_+-u^{n+1}_-), \\
\frac{\Delta x}{30}\dd{}{t}
(2h^n + 16h^{n+1/2} + 2h^{n+1}) & =& \frac{1}{6}(4u^n_+-4u^{n+1}_-),
\label{semi half}
\end{eqnarray}
\normalsize
where $h^n$ is the value of $h$ at the grid point $x^n$, $h^{n+1/2}$
is the value of $h$ at the midpoint $x^{n+1/2}$, $u^n_-$ is the
discontinuous $u$ value to the left of $x^n$ and $u^n_+$ is the value
to the right.

\begin{figure}
\begin{center}
\includegraphics*[width=14cm]{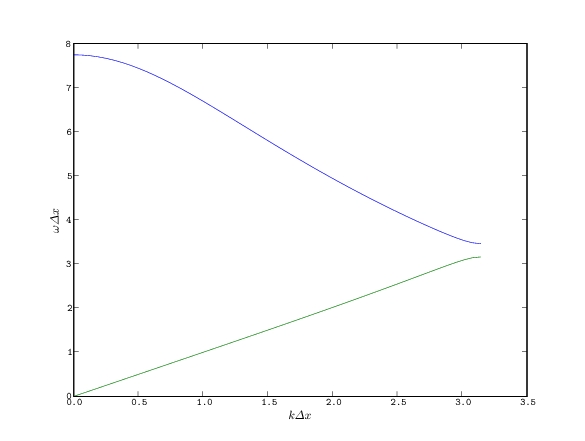}
\caption{\label{1d_eigs} Plot of the dispersion relation for the
semi-discrete equations obtained from the \pdgp element in one
dimension. The eigenspectrum has two branches, with a spectral gap
separating small and large eigenvalues. The lower branch is very straight
(and hence accurate).}
\end{center}
\end{figure}

We can obtain a dispersion relation for the semi-discrete system
(\ref{semi u+}-\ref{semi half}) by substituting the \emph{ansatz}
\begin{eqnarray*}
u^n_+ = \hat{u}_+e^{i(kx^n-\omega t)}, & & u^n_- = \hat{u}_-e^{i(kx^n-\omega t)}, 
\\
h^n = \hat{h}e^{ik(x^n-\omega t)}, 
& & h^{n+1/2} = \tilde{h}e^{i(kx^{n+1/2}-\omega t)}.
\end{eqnarray*}
We obtain the matrix equation
\begin{equation}
\small
\label{matrix eqn}
 \left[ \begin {array}{cccc} -2\,iw&-iw{e^{i\phi}}&-5+{e^{i\phi}}&4
\,{e^{1/2\,i\phi}}\\\noalign{\medskip}-iw&-2\,iw{e^{i\phi}}&-1+5\,{
e^{i\phi}}&-4\,{e^{1/2\,i\phi}}\\\noalign{\medskip}25-5\,{e^{-i\phi}}&
-25+5\,{e^{i\phi}}&-iw \left( 8-2\cos\phi\right) &-4
iw\cos{\phi/2}
\\\noalign{\medskip}-20&20\,{e^{i\phi}}&-2\,iw \left( 1+{e^{i\phi}}
 \right) &-16\,iw{e^{1/2\,i\phi}}\end {array} \right] 
\begin{pmatrix}
\hat{u}^+ \\
\hat{u}^- \\
\hat{h} \\
\tilde{h} \\
\end{pmatrix}
=
\begin{pmatrix}
0 \\ 
0 \\
0 \\
0 \\
\end{pmatrix}
\normalsize
\end{equation}
where $\phi = k\Delta x$ and $w=\omega\Delta x$.
After some algebraic manipulation using Maple, this yields
\[
w = \pm 2\sqrt{\frac{26 + 4\cos(\phi) 
\pm\sqrt{474 + 448\cos(\phi) - 22\cos(2\phi)}}
{6-2\cos\phi}}.
\]
A plot of this numerical dispersion relation is given in figure
\ref{1d_eigs}. The eigenvalues in the lower branch are monotonically
increasing, and there is a gap in the spectrum at $k\Delta x=\pi$. The
eigenvalues do not return to zero in the upper branch. The numerical
dispersion relation indicates that there are no spurious modes in the
discretisation and so the element is stable. Another feature is that
the low frequency branch is very close to the exact dispersion
relation for the wave equation.

\begin{figure}
\begin{center}
\includegraphics*[width=10cm]{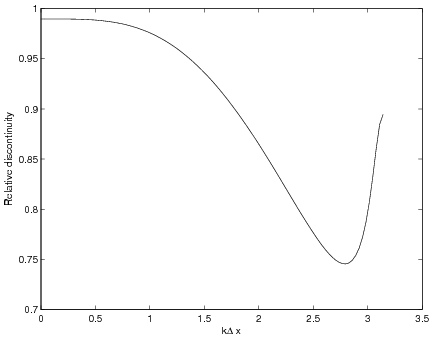}
\includegraphics*[width=10cm]{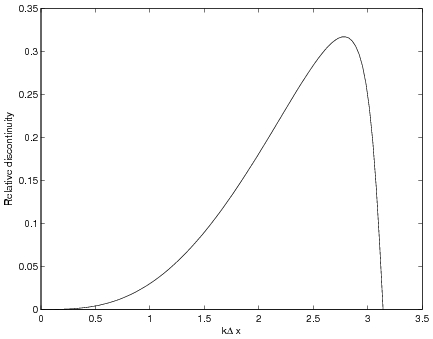}
\caption{\label{disc}Plots of the magnitude of the discontinuity in $u$
of the eigenmodes for the low (bottom plot) and high (top plot)
frequency branches of the dispersion relation. The low frequency modes
exhibit low levels of relative discontinuity and the high frequency
modes are very discontinuous with the fastest mode being completely
out of phase.}
\end{center}
\end{figure}

To investigate this gap in the spectrum further, we used this solution
to recover the structure of the modes by looking at the eigenvectors
of the matrix in equation (\ref{matrix eqn}) when $\omega$ takes these
values. We normalised the eigenvectors and calculated the magnitude of
the difference between $\hat{u}^+$ and $\hat{u}^-$, which gives a
measure of the discontinuity in each mode.  Figure \ref{disc}
illustrates that the level of discontinuity for modes from the lower
frequency branch is much lower than for those from the higher
frequency branch.

\section{Analysis of discrete Laplacian for two and three dimensional
  unstructured meshes }
\label{laplace}
In this section we construct the discrete Laplacian using the \pdgp
element for unstructured meshes in two and three dimensions and check the
eigenvalues for spurious modes. The meshes were produced by
Triangle\cite{Shewchuk96}
 (in two dimensions) and Tetgen\cite{Si2005} (in
three dimensions). The matrices $C$, $M^u$ and $M^h$ were assembled using
exact quadrature and the eigenvalues of the discrete Laplacian
$(M^h)^{-1}\cmc$ were computed numerically using Scientific Python.

\subsection{Two dimensions}
Table \ref{2D eigs Dirichlet} shows the computed eigenvalues for the discrete
Laplacian obtained from the \pdgp element in two dimensions with 
Neumann boundary conditions for $\vec{u}$ and Dirichlet boundary 
conditions for $h$. The meshes are unstructured in a $1\times1$ 
square domain.
\begin{table}
\begin{center}
\begin{tabular}{|c|c|c|}
\hline 
\begin{tabular}{c}
Max \\triangle area \\
\end{tabular}
& 
\begin{tabular}{c}
Number \\
of elements \\
\end{tabular}
& 
\begin{tabular}{c}
vector of eigenvalues \\
in
increasing magnitude \\
\end{tabular}
\\ 
\hline 
0.1 & 14 & $\begin{array}{c}
[19.44,40.79,44.95,55.42, \\
\qquad 61.42,...,583.7,656.3]\end{array}$ \\ 
\hline 
0.05 & 36 & $\begin{array}{c}
[19.72,48.98,49.16,76.34, \\
\qquad 94.41,...,9102,29160] \\
\end{array}$ \\ 
\hline 
0.01 & 151 & $\begin{array}{c}
[19.73,49.35,49.36,79.02, \\
\qquad 98.80,...126300,2949000] \\
\end{array}$ \\
\hline 
\end{tabular}
\end{center}
\caption{\label{2D eigs Dirichlet} Table showing eigenvalues of the
discrete Laplacian obtained from the \pdgp element pair in two
dimensions with Dirichlet boundary conditions. All the eigenvalues
correspond to physical modes.}
\end{table}

The Dirichlet boundary conditions for $h$ prohibit the constant $h$
solution with eigenvalue zero and so the smallest physical eigenvalue
is $2\pi^2$ corresponding to $h=\sin(x)\sin(y)$. It is clear from the
table that there are no spurious eigenvalues \emph{i.e.} eigenvalues
that scale with the mesh size, and all of the eigenvalues correspond
to physical modes.  

Table \ref{2D eigs Neumann} shows the computed
eigenvalues for the discrete Laplacian obtained from the \pdgp element
in two dimensions with Neumann boundary conditions for $h$ and
Dirichlet boundary conditions for $\vec{u}$, on the same $1\times1$
domain.
\begin{table}
\begin{center}
\begin{tabular}{|c|c|c|}
\hline 
\begin{tabular}{c}
Max \\triangle area \\
\end{tabular}
& 
\begin{tabular}{c}
Number \\
of elements \\
\end{tabular}
& 
\begin{tabular}{c}
vector of eigenvalues \\
in
increasing magnitude \\
\end{tabular}
\\ 
\hline 
0.1 & 14 & 
$\begin{array}{c}
[0.00,9.89,9.90,19.90, \\
\qquad 41.22,...,830.0,960.9]\\
\end{array}$ \\
\hline 
0.05 & 36 & 
$\begin{array}{c}[0.00,9.88,9.88,19.80, \\
\qquad 40.43,...,15870,41690] \\
\end{array}$\\
\hline 
0.01 & 151 & 
$\begin{array}{c}[0.00,9.87,9.87,19.74, \\
\qquad 39.50,...,177300,3274000] \\
\end{array}$ \\
\hline 
\end{tabular}
\end{center}
\caption{\label{2D eigs Neumann} Table showing eigenvalues of the
discrete Laplacian obtained from the \pdgp element pair in two
dimensions with Neumann boundary conditions. All the eigenvalues
correspond to physical modes.}
\end{table}
The Neumann boundary conditions for $h$ admit the constant $h$
solution with eigenvalue zero. The next two physical eigenfunctions
are $\sin(\pi x)$ and $\sin(\pi y)$ which both have eigenvalues
$\pi^2$. There are no spurious eigenvalues.

\subsection{Three dimensions}
Table \ref{3D eigs Dirichlet} shows the computed eigenvalues for the discrete
Laplacian obtained from the \pdgp element in three dimensions with 
Neumann boundary conditions for $\vec{u}$ and Dirichlet boundary 
conditions for $h$. The meshes are unstructured in a $1\times1\times1$ 
cubic domain.
\begin{table}
\begin{center}
\begin{tabular}{|c|c|c|}
\hline 
\begin{tabular}{c}
Max \\tetrahedral volume \\
\end{tabular}
& 
\begin{tabular}{c}
Number \\
of elements \\
\end{tabular}
& 
\begin{tabular}{c}
vector of eigenvalues \\
in
increasing magnitude \\
\end{tabular} \\
\hline 
0.1 & 44 & 
$\begin{array}{c}
[0.00,0.00,0.00,0.00,0.00, \\
\qquad 29.876,...,939.7,1045] \\
\end{array}$
\\
\hline 
0.01 & 215 & $\begin{array}{c}
[0.00,0.00,0.00,0.00,0.00, \\
\qquad 0.00,29.87,...,5177,5753] \\
\end{array}$\\
\hline 
0.0059 & 329 & 
$\begin{array}{c}
[0.00,0.00,0.00,29.71, \\
\qquad 59.85,...,5772,6723]
\end{array}$ \\
\hline 
0.00585 & 330 & 
$\begin{array}{c}
[29.72,59.88,60.06, \\
\qquad 60.07,...,5763,6631] \\
\end{array}$ \\
\hline
0.005 & 398 & 
$\begin{array}{c}
[29.72,59.88,60.06, \\
\qquad 60.07,...,7474,11570] \\
\end{array}$\\
\hline
0.004 & 487 & 
$\begin{array}{c}
[29.70,59.80,59.87,59.96, \\
\qquad 90.77,...,10680,12800] \\
\end{array}$\\
\hline
0.003 & 681 & 
$\begin{array}{c}
[29.67,59.57,59.61,59.66, \\
\qquad 90.15,...,11610,12350] \\
\end{array}$\\
\hline
\end{tabular}
\end{center}
\caption{\label{3D eigs Dirichlet} Table showing eigenvalues of the
discrete Laplacian obtained from the \pdgp element pair in three
dimensions with Dirichlet boundary conditions. There are spurious
eigenvalues for the coarsest meshes which disappear when there are more
elements and the ratio of $\vec{u}$ DOF to $h$ DOF is greater.}
\end{table}

As in the two-dimensional case, the Dirichlet boundary conditions for
$h$ prohibit the constant $h$ solution with eigenvalue zero and so the
smallest physical eigenvalue is $3\pi^2$ corresponding to
$h=\sin(x)\sin(y)\sin(z)$. Table \ref{3D eigs Dirichlet} shows that 
spurious eigenvalues are present for very coarse meshes but are not present
when the number of elements is increased.

Table \ref{3D eigs Neumann} shows the computed eigenvalues for the discrete
Laplacian obtained from the \pdgp element in three dimensions with
Neumann boundary conditions for $h$ and Dirichlet boundary conditions
for $\vec{u}$. The meshes are unstructured in a $1\times1\times1$
cubic domain.
\begin{table}
\begin{center}
\begin{tabular}{|c|c|c|}
\hline 
\begin{tabular}{c}
Max \\tetrahedral volume \\
\end{tabular}
& 
\begin{tabular}{c}
Number \\
of elements \\
\end{tabular}
& 
\begin{tabular}{c}
vector of eigenvalues \\
in
increasing magnitude \\
\end{tabular} \\
\hline 
0.1 & 44 & 
$\begin{array}{c}
[0.00,9.93,9.93,10.06, \\
\qquad 20.15,...,984,1097] \\
\end{array}$ \\
\hline 
0.01 & 215 & 
$\begin{array}{c}
[0.00,9.88,9.88,9.89, \\
\qquad 19.83,...,5385,5931] \\
\end{array}$\\
\hline
0.005 & 398 & 
$\begin{array}{c}
[0.00,9.874,9.874,9.875, \\
\qquad 19.78,...,7746,12070] \\
\end{array}$ \\
\hline
0.004 & 487 & 
$\begin{array}{c}
[0.00,9.873,9.873,9.873, \\
\qquad 19.78,...,10780,13010] \\
\end{array}$\\
\hline 
0.003 & 681 & 
$\begin{array}{c}
[0.00,9.872,9.872,9.873, \\
\qquad 19.76,...,11600,12330] \\
\end{array}$\\
\hline 
\end{tabular}
\end{center}
\caption{\label{3D eigs Neumann} Table showing eigenvalues of the
discrete Laplacian obtained from the \pdgp element pair in three
dimensions with Neumann boundary conditions. All the eigenvalues
correspond to physical modes, indicating that the element pair is
stable.}
\end{table}

The Neumann boundary conditions for $h$ admit the constant $h$
solution with eigenvalue zero. The next three physical eigenfunctions
are $\sin(\pi x)$, $\sin(\pi y)$ and $\sin(\pi z)$ which both have
eigenvalues $\pi^2$. There are no spurious eigenvalues.

\section{Numerical test for the wave equation}
\label{wave eqn}
In this section we test the \pdgp element as applied to the 
wave equation in two dimensions, with the aim of checking that spurious
oscillations do not appear and that the solution remains smooth.

We discretised the equations in time using the St\"ormer-Verlet method
given by
\begin{eqnarray*}
M^u\frac{\U_i^{n+1/2}-\U_i^n}{2\Delta t} & = &
-C_i\h^n, i=1,\ldots d, \\
M^h\frac{\h^{n+1}-h^n}{\Delta t} & = & 
\sum_{i=1}^dC^T_i\U^{n+1/2}, \\
M^u\frac{\U_i^{n+1}-\U_i^{n+1/2}}{2\Delta t} & = &
-C_i\h^{n+1}, i=1,\ldots d, 
\end{eqnarray*}
This method is second-order in time, and is symplectic, one of the
consequences of which is that there exists a conserved energy which is
equal to the exact spatially discretised energy plus a correction of
magnitude $\mathcal{O}(\Delta t^2)$ (see \cite{LeRe2005} for a review
of the St\"ormer-Verlet method applied to PDEs). This means that
small-scale energy will not be dissipated and it provides a good test
of the spatial discretisation. As this method is explicit, there is a
numerical CFL condition which requires that the fastest oscillation in
the system, corresponding to the largest eigenvalue of the discrete
Laplacian, should be resolved in time. This discretisation still
requires linear systems to be solved to obtain $\vec{u}$ and $h$ at
the next time level, although the mass matrix for $\vec{u}$ is block
diagonal (with one block per element).

\begin{figure}
\begin{center}
\includegraphics*[width=6cm]{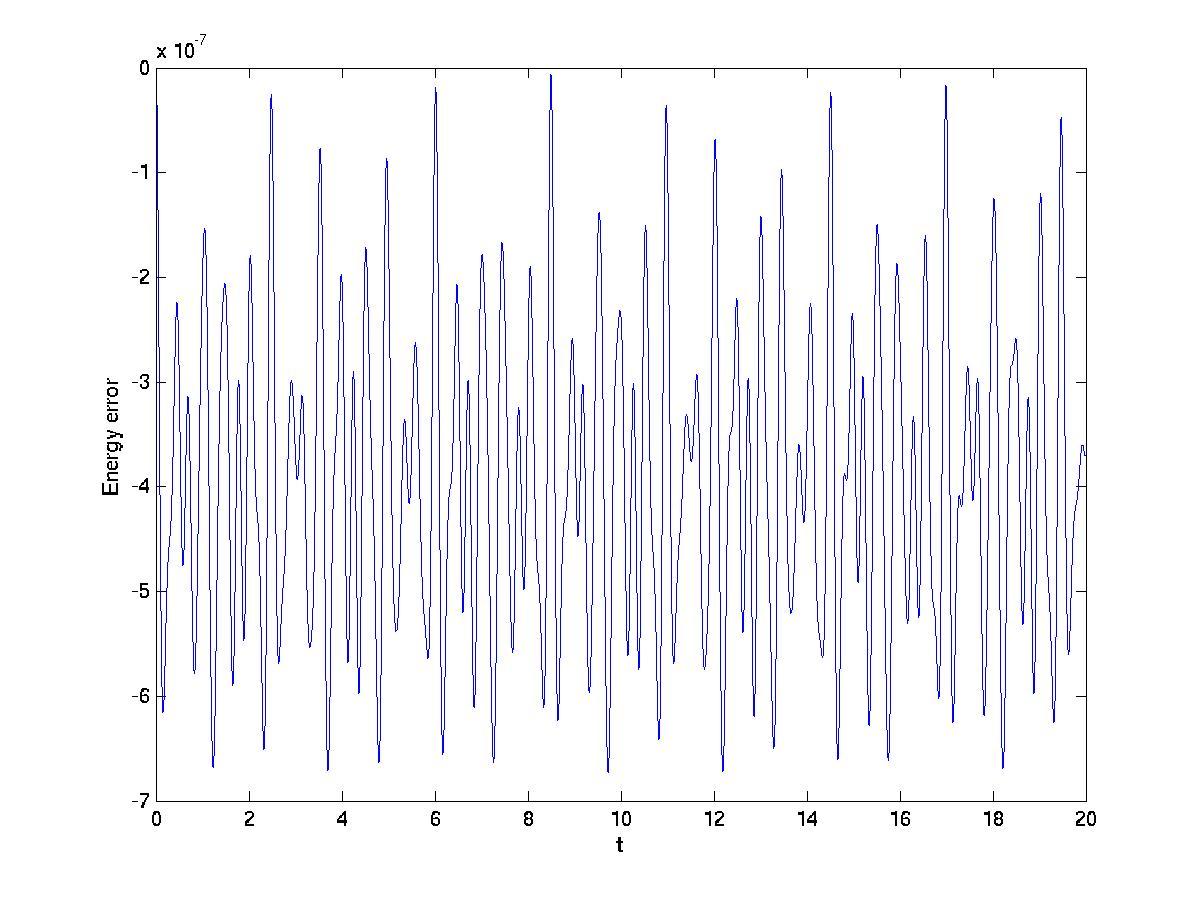}
\includegraphics*[width=6cm]{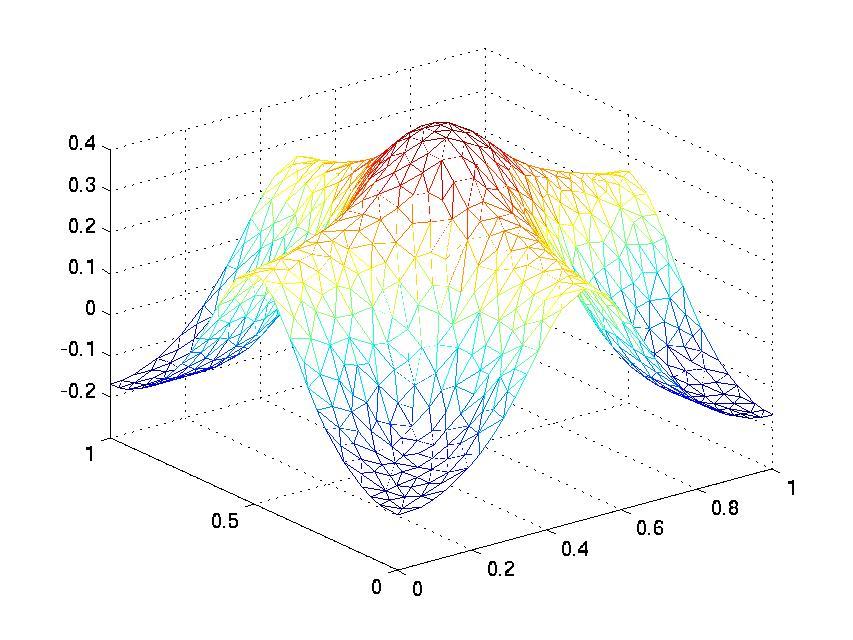} \\
\includegraphics*[width=6cm]{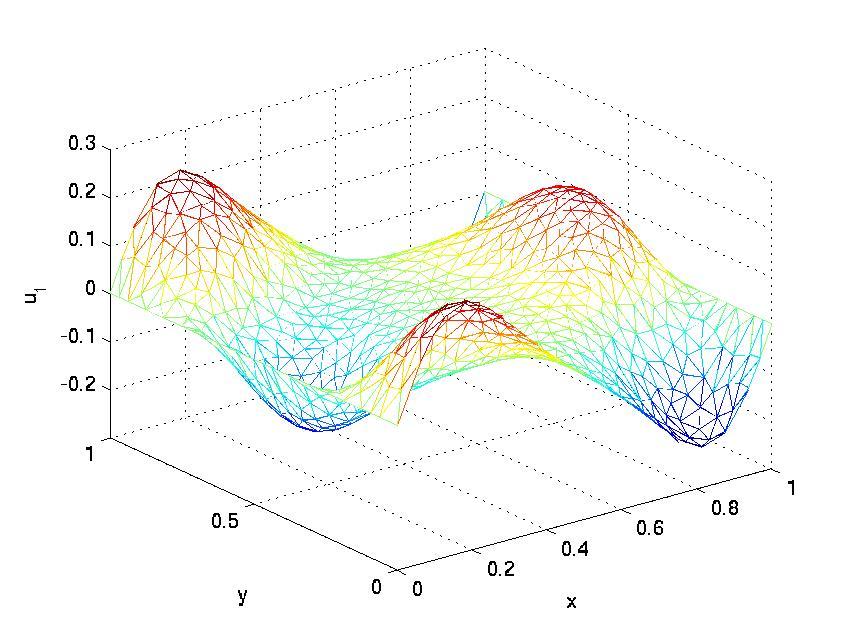} 
\includegraphics*[width=6cm]{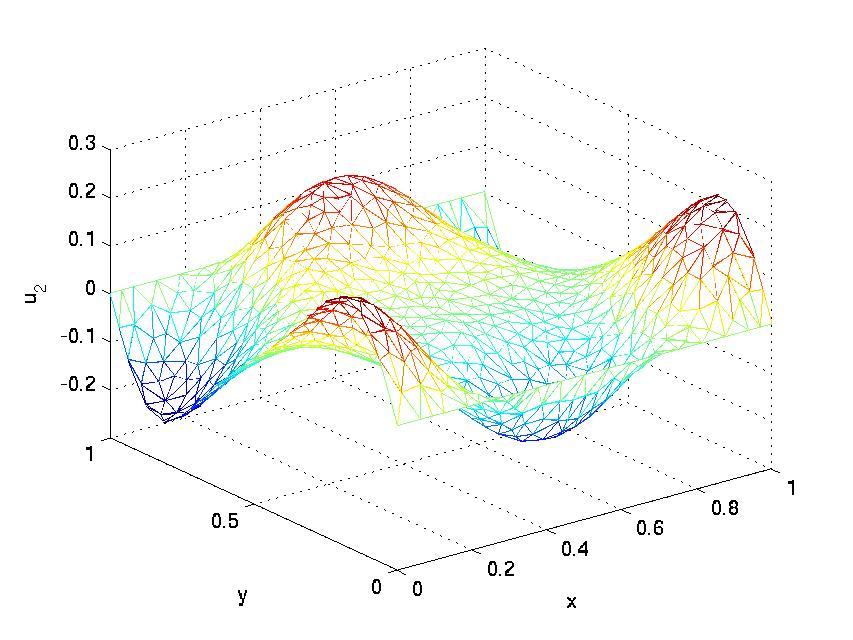}
\end{center}
\caption{\label{wave fig}Numerical results obtained from solving the
  2-dimensional wave equation in a square domain of area 1 with an
  unstructured grid with triangular elements of typical area
  0.001. The wave speed is $c=1$ and the timestep is $\Delta
  t=0.001$. Top-left: plot of energy error against time. Top-right:
  plot of $h$ at time $t=20.0$. Bottom-left and bottom-right: plots of
  the $x$- and $y$-components of $\vec{u}$. These results show that
  the numerical solution stays smooth after a large number of
  timesteps, which is a good indicator that the method is
  stable. Although it appears from the plot that $\vec{u}$ 
  remains almost continuous, small discontinuities are present
  in the solution.}

\end{figure}

Simulation results are given in figure \ref{wave fig}. These results
show that the solutions remain smooth and that there are no spurious modes
polluting the solution. This good behaviour arises from the stability of
the \pdgp element.

\section{Summary and Outlook}
\label{summary}
In this paper we introduced the \pdgp mixed element which has
discontinuous velocity and continuous pressure. This choice means that
the number of DOF for velocity can be increased in order to obtain a
stable element. In section \ref{mixed} we described the construction
of the element in one, two and three dimensions and gave example
values for the $\vec{u}$ and $h$ DOF. In future implementations in the
three-dimensional non-hydrostatic Imperial College Ocean Model (ICOM)
\cite{PaPiGo2005} we will investigate the relative merits of \pdgp,
\pdg-$\textrm{P3}$, $\textrm{P2}_{\textrm{DG}}$-$\textrm{P3}$ and
other combinations, including augmenting the $\vec{u}$ space with
bubble functions, in practical applications.

In section \ref{1d} we gave a linear normal mode analysis for the
element on a regular grid in one dimension with periodic boundary
conditions which showed that the element is stable in this case. The
dispersion relation is monotonically increasing with a spectral gap
between the two branches, and the lower frequency branch has
relatively continuous eigenfunctions with almost continuous
eigenfunctions at the lowest frequencies.

In section \ref{laplace} we presented calculations of eigenvalues of
the discrete Laplace matrix obtained from unstructured meshes in two
and three dimensions which demonstrated that the element is stable in these
cases. In section \ref{wave eqn} we presented results from a wave equation
calculation on a two dimensional grid which demonstrated that the solutions 
stay smooth for relatively long time intervals.

This type of element with discontinuous velocity and continuous pressure has
some other properties that may make it favourable for use in geophysical
codes such as ICOM. The discontinuous element for velocity means that the
discretisation locally conserves momentum, and the use of upwinding and
flux-limiting means that the discretisation allows a good treatment of
advection (see, for example, \cite{Levin2006}). As the mass matrix for
$\vec{u}$ is block diagonal, the $\cmc$ matrix remains sparse and so it is
not necessary to lump the mass matrix. This makes the discretisation more
accurate and reduces problems with balancing various lumped and non-lumped
terms.

A key issue with modelling large-scale rotating geophysical flow is that
of geostrophic balance, which states that the Coriolis term is nearly
balanced by the horizontal pressure gradients:
\[
\vec{\Omega}\times\vec{u}\approx -\nabla_Hp,
\]
where $\vec{\Omega}$ is the Earth's rotation vector and $\nabla_H$ is
the horizontal gradient. For an element pair such as P1-P1, the
pressure gradients are piecewise constant whilst the Coriolis force is
piecewise linear and it is not possible to find a pressure field to
accurately represent this balance. This leads to pressure gradient
errors which pollute the solution after short times, and it becomes
necessary to subtract out the balanced pressure (discretised on a
higher-order element) in order for the solution to stay near to
balance. For the \pdgp element pair, the velocity field is piecewise
linear whilst the pressure field is piecewise quadratic, and it will
be possible to find a pressure field which represents this balance as
long as the velocity field remains relatively continuous. The study of
discontinuity in the normal modes for the one-dimensional element is
encouraging as it shows that the lower branch of the spectrum,
corresponding to well-resolved solutions, remains relatively
continuous. We will investigate the pressure gradient errors arising
from this discretisation in future work.

\nocite{*}
\bibliography{DgCg}

\end{document}